\documentclass[9pt]{article}
\usepackage{amsmath}
\usepackage{latexsym}
\usepackage{amsfonts}
\usepackage{amssymb}
\usepackage{makeidx}

\newcommand{\Ss}{\cS_s}
\newcommand{\Sf}{\cS_f}
\newcommand{\Ls}{\cL_s}
\newcommand{\Lf}{\cL_f}
\newcommand{\Qf}{\cQ_f}
\newcommand{\Lc}{\cL_f^c}
\newcommand{\Qc}{\cQ_f^c}

\def\la{\langle}
\def\ra{\rangle}

\def\ll{\langle\langle}
\def\rr{\rangle\rangle}

\renewcommand{\phi}{\varphi}

\newcommand{\nf}{\hspace*{-5pt}}

\def \cao {\c c\~ao}

\def \beq { \begin{equation} }
\def \eeq { \end{equation} }

\def\overset#1\to#2{\mathrel{\mathop{#2}\limits^{#1}}}
\def\underset#1\to#2{\mathrel{\mathop{#2}\limits_{#1}}}

\def \varinjlim {\underset{\longrightarrow}\to{lim}}
\def \varprojlim {\underset{\longleftarrow}\to{lim}}

\def \li {\varinjlim}
\def \lp {\varprojlim}

\def \lli#1 {\mathrel{\mathop{\li}\limits_{#1}}}
\def \llp#1 {\mathrel{\mathop{\lp}\limits_{#1}}}

\def \rest {\restriction}
\def \la {\leftarrow}
\def \im {\Rightarrow}

\def \rsa {\rightsquigarrow }
\def \lrsa {\leftrightsquigarrow }

\def \thra {\twoheadrightarrow }

\def \hookr {\hookrightarrow }

\def \rlas {\rightleftarrows}

\newcommand{\hem}{\hspace*{1em}}

\newcommand{\hfl}{\hspace*{\fill}}

\newlength{\dede}     

\newcommand{\hp}{\hspace*{\parindent}} 

\newcommand{\N}{\mbox{$ \mathbb N$}}   

\newcommand{\lra}{\mbox{$\longrightarrow$}}  

\newcommand{\sub}{\mbox{$\subseteq$}}

\newcommand{\cK}{\mbox{$\cal K$}}
\newcommand{\cL}{\mbox{$\cal L$}}

\newcommand{\cA}{\mbox{$\cal A$}}

\newcommand{\cQ}{\mbox{$\cal Q$}}

\newcommand{\cS}{\mbox{$\cal S$}}

\newcommand{\cX}{\mbox{$\cal X$}}
\newcommand{\cY}{\mbox{$\cal Y$}}

\newcommand{\all}{\mbox{$\forall$}}

\newcommand{\Ra}{\mbox{$\Rightarrow$}}
\newcommand{\Lra}{\mbox{$\Leftrightarrow$}}

\def \rest {\restriction}
\def \ra {\rightarrow}
\def \la {\leftarrow}
\def \im {\Rightarrow}

\def \rsa {\rightsquigarrow }
\def \lrsa {\leftrightsquigarrow }

\def \thra {\twoheadrightarrow }

\def \hookr {\hookrightarrow }

\def\overset#1\to#2{\mathrel{\mathop{#2}\limits^{#1}}}
\def\underset#1\to#2{\mathrel{\mathop{#2}\limits_{#1}}}

\newcommand{\und}[1]{\raisebox{-.2ex}{\underline{\raisebox{.2ex}{#1}}}}  

\def \ll {\langle}
\def \rr {\rangle}

\def \varinjlim {\underset{\longrightarrow}\to{lim}}
\def \varprojlim {\underset{\longleftarrow}\to{lim}}

\def \li {\varinjlim}
\def \lp {\varprojlim}

\def \lli#1 {\mathrel{\mathop{\li}\limits_{#1}}}
\def \llp#1 {\mathrel{\mathop{\lp}\limits_{#1}}}

\newtheorem{Th}{Theorem}[section] 
\newtheorem{Co}[Th]{Corollary}
\newtheorem{Df}[Th]{Definition}
\newtheorem{Pro}[Th]{Proposition}
\newtheorem{Le}[Th]{Lemma}
\newtheorem{Exa}[Th]{Example}
\newtheorem{Rem}[Th]{Remark}
\newtheorem{Fa}[Th]{Fact}
\newtheorem{Que}[Th]{Question}
\newtheorem{Ct}[Th]{}
\newtheorem{Afi}[Th]{Claim}

\newcommand{\baf}{\begin{Afi}\nf{\sl . }}
\newcommand{\eaf}{\end{Afi}}

\newcommand{\bdf}{\begin{Df}\nf{\bf .}}
\newcommand{\edf}{\end{Df}}
\newcommand{\bte}{\begin{Th}\nf{\bf .}}
\newcommand{\ete}{\end{Th}}  
\newcommand{\bco}{\begin{Co}\nf{\bf .}}
\newcommand{\eco}{\end{Co}}
\newcommand{\ble}{\begin{Le}\nf{\bf .}}
\newcommand{\ele}{\end{Le}}
\newcommand{\bpr}{\begin{Pro}\nf{\bf .}}
\newcommand{\epr}{\end{Pro}}
\newcommand{\bex}{\begin{Exa}\nf{\bf .} \rm}
\newcommand{\eex}{\end{Exa}}
\newcommand{\bre}{\begin{Rem}\nf{\bf .} \rm}
\newcommand{\ere}{\end{Rem}}
\newcommand{\bfa}{\begin{Fa}\nf{\bf .} \sl}
\newcommand{\efa}{\end{Fa}}
\newcommand{\bqt}{\begin{Que}\nf{\bf .}}
\newcommand{\eqt}{\end{Que}}

\newcommand{\qdr}{\hfl $\square$ }
\newcommand{\bct}{\begin{Ct}\nf{\bf .} \rm}
\newcommand{\ect}{\end{Ct}}
\newcommand{\bxa}{\begin{Exa}\nf{\bf .} \rm}
\newcommand{\exa}{\end{Exa}}

\newcommand{\bu}{$\bullet$\ }

\begin{document}

\title{Representation theory of logics: a categorial approach} 

\author{Darllan Concei\cao\ Pinto \thanks{Instituto de Matem\'atica e Estat\'istica, University of S\~ao Paulo, Brazil.\ Emails: darllan\_\! math@hotmail.com, hugomar@ime.usp.br }\\
Hugo Luiz Mariano \thanks{Research supported by FAPESP, under the Thematic Project LOGCONS:
Logical consequence, reasoning and computation (number 2010/51038-0).} 
}

\date{December 2013}

\maketitle

\begin{abstract}
The major concern in the study of categories of logics  is to describe condition for preservation, under the a method of combination of logics, of meta-logical  properties. Our complementary approach to this field is study the "global" aspects of categories of logics in the vein of the categories $\Ss, \Ls, \cA_s$ studied in \cite{AFLM3}. All these categories have good properties however the category of logics $\cL$ does not allow a good treatment of the "identity problem" for logics (\cite{Bez}):   for instance, the presentations of "classical logics" (e.g., in the signature $\{\neg, \vee\}$ and $\{\neg',\rightarrow'\}$) are not $\Ls$-isomorphic. In this work, we  sketch  a possible way to overcome this "defect" (and anothers) by a mathematical device: a representation theory of logics obtained from category theoretic aspects on (Blok-Pigozzi) algebraizable logics. 
In this setting we propose the study of (left and right) "Morita equivalence" of logics and variants.
We introduce the concepts of logics (left/right)-(stably) -Morita-equivalent and show  that the presentations of classical logics are stably Morita equivalent but classical logics and intuitionist logics are not stably-Morita-equivalent: they are only stably-Morita-adjointly related. 
\end{abstract}

\section{Introduction}

\hp In the 1990's rise many methods of combinations of logics (\cite{CC3}). They appear in dual aspects: as processes of decomposition or analysis of logics (e.g., the "Possible Translation Semantics" of W. Carnielli, , \cite{Car}) or as a processes  of  composition or synthesis of logics (e.g., the "Fibrings" of D. Gabbay, \cite{Ga}). This was the main motivation  of categories of logics. The major concern in the study of categories of logics (CLE-UNICAMP, IST-Lisboa) is to describe condition for preservation, under the combination method, of meta-logical  properties (\cite{CCCSS}, \cite{ZSS}). Our complementary approach to this field is study the "global" aspects of categories of logics (\cite{AFLM1}, \cite{AFLM2}, \cite{AFLM3}, \cite{MaMe}).

The initial steps on "global" approach to categories of logics are given in the sequence of papers \cite{AFLM1}, \cite{AFLM2} and \cite{AFLM3}: they present very simple but too strict notions of  logics and morphisms, with "good" categorial properties (\cite{AR}) but unsatisfactory treatment of the "identity problem" of logics (\cite{Bez}). More  flexible notions of morphisms between logics are considered in \cite{FC}, \cite{BCC1}, \cite{BCC2}, \cite{CG}: this alternative notion allows better approach to the identity  problem however has many categorial "defects". A "refinement" of those ideas is provided in  \cite{MaMe}: are considered categories of logics satisfying {\em simultaneously} certain natural  conditions: (i) represent the major part of logical systems; (ii) have  good categorial properties; (iii) allow  a natural notion of  algebraizable  logical system (\cite{BP}, \cite{Cze}); (iv) allow satisfactory treatment of the "identity problem" of logics.

In the present work we present and alternative approach to overcome the  problems above through a {\em mathematical device}:  a representation theory of logics obtained from category theoretic aspects on (Blok-Pigozzi) algebraizable logics:

{\bf Motivation 1}: analogy: logics $\lrsa$ rings\\
\bu "Representation theory of rings":\\
-- $R \in obj(Ring) \ \rsa$ \ $R-Mod$ (respec., $Mod-R$) $\in CAT$;\\
-- (left/right) Morita equivalence of rings: \\ $R \equiv R'$ \ $\Lra$ \ 
$R-Mod \simeq R'-Mod$  (respec., $Mod-R \simeq Mod-R'$).\\
\bu "Representation theory of propositional logics":\\
-- $l \in obj(Log) \ \rsa$ \ $l-Mod$ (respec., $Mod-l$): diagrams of categories and functors (respec.: diagrams of categories, functors and natural transformations);\\
-- (left/right) Morita equivalence of  logics and variants;\\
-- left and right are conceptually and technically distincts.

{\bf Motivation 2}:  analogy: logics $\lrsa$ topology\\
\bu "Algebraic Topology": (objects: topological spaces or logics) \\
-- define a general theory of "mathematical invariants" to mesure the degree of distinctions of arbitrary logics;\\
-- develope general methods of calculation of "invariants" (in some sense); \\
-- introduce new forms of comparation of objects.

\section{Preliminaries}

\subsection{Categories of Signatures and Categories of Logics}

\hp If we want define and study  categories of logics, we must provide answers to the two natural questions: (i) how to represent a logical system? (ii) what are the relevant notions de morphisms? (\cite{CC1}, \cite{CC2}, \cite{CCRS}, \cite{SSC}). In the following  we adopt a simple --and sintatical-- approach to this theme.

\bct \ $\Ss$, the category of signature and (strict or simple) signature morphisms:\\
A  (propositional, finitary) signature is a sequence of pairwise disjoint sets $\Sigma = (\Sigma_{n})_{n \in \N}$. In what follows,  $X = \{x_{0},x_{1},\ldots,x_{n},\ldots\}$
will denote a fixed enumerable set (written in a fixed order). Denote  $F(\Sigma)$ (repectively $F(\Sigma)[n]$), the set of $\Sigma$-formulas over build from $X$ (respec. $\{x_0, \ldots, x_{n-1}\}$). A (strict) morphism, $f : \Sigma \ra \Sigma'$, is a sequence of functions:
$f = (f_n)_{n \in \N} : (\Sigma_{n})_{n \in \N} \ra (\Sigma'_{n})_{n \in \N}$; this induces a function between formula algebras  $\hat{f} : F(\Sigma) \ra F(\Sigma')$.
\qdr\ect

\bct \ $\Ss \simeq Set^{\N}$, is a finitely locally presentable  category  and the fp signatures are the "finite support" signatures.
\qdr\ect

Recall:\\
{{\em \small {(i) locally presentable (= accessible + complete/cocomplete) (\cite{AR}, \cite{MP});\\
(ii) a category is $\kappa$-accessible if it has $\kappa$-filtered colimits and a {\em set} of $\kappa$-presentable objects such that every object is a $\kappa$-filtered colimit of these objects}}}\

\bct \ $\Ls$, the category of (strict) logics over $\Ss$:\\
A logic is a pair $l = (\Sigma, \vdash)$, where $\Sigma$ is a signature and 
$\vdash$ is a tarskian consequence operator.  A $\Ls$-morphism,  $f : (\Sigma,\vdash) \ra (\Sigma', \vdash')$, is a (strict) signature morphism 
$f \in \Ss(\Sigma,\Sigma')$ such that
$\hat{f} : F(\Sigma) \ra F(\Sigma')$ is a $(\vdash,\vdash')$-translation:  $\Gamma \vdash \psi \ \Ra \ \hat{f}[\Gamma] \vdash' \hat{f}(\psi)$, for all $\Gamma \cup\{\psi\} \sub F(\Sigma)$ (i.e., it is "continuous"). 
\qdr\ect


\bct \ $\Ls$ is a $\omega$-locally  presentable category and the fp logics are given by a finite set of  "axioms" and "inference rules" over a fp signature.
\qdr\ect



\bct \ $\cA_s$, the category (strict) of BP-algebraizable logics (see \cite{BP}):\\
\bu objects: logic $l = (\Sigma, \vdash)$, that has some {\em algebraizing pair} $((\delta\equiv\epsilon), \Delta)$;\\
\bu morphisms:  $f: l \ra l'$ : $f \in \Ls(l,l')$  and   ``preserves algebraizing  pair'' (well defined).
\qdr\ect

 Recall that:\\
 {\small {   $\{\delta_{r}, \epsilon_{r} : r  < s\} \ \sub \ F(\Sigma)[1]$;\\
    $\{\Delta_{u}: u<v\} \sub F(\Sigma)[2]$;\\
    $((\delta\equiv\epsilon), \Delta)$ satisfies conditions (i) and (ii) (and/or
conditions (i)' and (ii)') below, with $\Gamma \cup \Theta \cup \{\psi, \varphi,
\zeta, \eta, \vartheta\} \sub F(\Sigma)$:\\
    (i) $\Gamma\vdash \varphi \Lra \{(\delta(\psi)\equiv\epsilon(\psi)) :
    \psi \in \Gamma\} \vDash_{K} (\delta(\varphi)\equiv\epsilon(\varphi))$;\\
    (ii) $(\varphi \equiv \psi)$ $\vDash_{K}\vDash$
    $(\delta(\varphi\Delta\psi) \equiv
    \epsilon(\varphi\Delta\psi))$;\\
    (i)' $\Theta\vDash_{K} (\varphi\equiv\psi) \Lra \{\zeta\Delta\eta :
    (\zeta\Delta\eta) \in \Theta\} \vdash \varphi\Delta\psi$;\\
    (ii)' $\vartheta \dashv \vdash
    \delta(\vartheta)\Delta\epsilon(\vartheta)$.
}}

\bct \ Functors:\\
\bu  Forgetful functors: $U : \cA \ra \Ls$; \hem $U' : \cA \ra \Ss$;\\
\bu If $f \in \cA(l_0,l_1)$ and $\cK^i \sub \Sigma^i-Str$ is the quasivariety equivalent algebraic semantic of $l_i$, then $f^\star : \Sigma^1-Str \ra \Sigma^0-Str$ ($M_1 \mapsto (M_1)^f$) restricts to $f^\star : \cK^1 \ra \cK^0$.  
\qdr\ect

\bct \ Limits and colimits in $\cA_s$:\\
 $U$ creates products over  {\em ''bounded''} diagrams. $U'$ creates colimits over   {\em non empty} diagrams. 
$U$ creates {\em filtered colimits}, moreover 
if $(l, (\gamma_i)_{i \in (I,\leq)})$ is a colimit cocone, then given $M \in \Sigma-Str$, 
$M \in \cK$ $\Lra$ $M^{\gamma_i} \in \cK_i$, $\all i \in I$.
\qdr\ect

\bct \ $\cA_s$ is a finitely accessible category (but not complete/cocomplete). Moreover $U : \cA \ra \Ls$ is a $\omega$-accessible functor.
\qdr\ect

\bct \ Remote algebrization revisited (LFIs):\\
\bu {\small $F_i: \cX_i \ra \cY$, $i =0,1$, accessible functors $\im$ $(F_0 \ra F_1)$ is an accessible category;}\\
\bu {\small accessible categories have a {\em small} weakly initial family;}\\
\ {\em Proposition:} For each $l \in Obj(\Ls)$, there is a {\em small} family of $\Ls$-morphisms $(\eta_i : l \ra U(l_i))_{i \in I}$ such that for each $l' \in Obj(\cA)$ and $f \in \Ls(l, U(l'))$, there are $i \in I$ and $f_i \in \cA(l_i,l')$ such that $U(f_i) \circ \eta_i = f$.\\
\ {\em Corollary:} A {\em weak} universal property of $\eta : l \ra \prod_{i \in I} U(l_i)$. \\
\ {\em Questions:} Describe conditions on $l$ such that:\\
-- $\{l_i: i \in I\} \ \sub \ (\cA)_{fp}$.\\
-- $\{l_i : i \in I\}$ be bounded and $U(\prod_{i \in I} l_i) \overset{\cong}\to\ra \prod_{i \in I} U(l_i)$.\\
-- we can replace $\exists \lrsa \exists!$.\\
Then the ingredients are "canonical" ($I \cong I'$, $l_i \cong l_{i'}$) and allows us to define  
"the algebraizable spectrum of the logic $l$"\\
{\small (analagogy with rings: $R \in Obj(cRing_1) \rsa (\alpha_P : R \ra Frac(R/P))_{P \in Spec(R)}$)}.
\qdr
\ect

\bct \ But they does not allow a good treatment of the "identity problem" for logics:   for instance, the presentations of "classical logics" (e.g., in the signature $\{\neg, \vee\}$ and $\{\neg',\rightarrow'\}$) are not $\Ls$-isomorphic. 
\qdr\ect

{\bf {In this work, we  sketch  a possible way to overcome this "defect", by a mathematical device.}}

\bct {\em Other categories of logics}

$\bullet$ \ $\Lf$: logical translations with "flexible" signature morphisms  \
$c_n \in \Sigma_n \mapsto \varphi'_n \in F(\Sigma')[n]$ (\cite{FC})

$\bullet$ \ $Q\Lf$: "quotient" category: $f \sim g$ iff\
$\check{f}(\varphi) \dashv'\vdash \check{g}(\varphi)$.\\
The logics $l$ and $l'$ are equipollent (\cite{CG}) iff $l$ and $l'$ are $Q\Lf$-isomorphic.

$\bullet$ \ $\Lc \sub \Lf$: "congruential" logics: \
$\varphi_0 \dashv \vdash \psi_0, \ldots, \varphi_{n-1} \dashv \vdash  \psi_{n-1}$ $\Ra$\
 $c_{n}(\varphi_0, \ldots, \varphi_{n-1}) \dashv \vdash c_{n}(\psi_0, \ldots, \psi_{n-1})$.\\
 The inclusion functor $\Lc \hookr \Lf$ has a left adjoint.

$\bullet$ \ $Lind(\cA_f) \sub \cA_f$: "Lindenbaum algebraizable" logics:\
$ \varphi \dashv  \vdash \psi$ $\Lra$
        $\vdash \varphi\Delta\psi$ (well defined).\\
$Lind(\cA_f) \sub \Lc$ and the inclusion functor $Lind(\cA_f) \hookr \cA_f$ has a left adjoint.

$\bullet$ \ $Q\Lc$ (or simply $\Qc$): "good" category of logics: represents the major part of logics; has good categorial properties (is an accessible category complete/cocopmplete); solves the identity problem for the presentations of classical logic interms of isomorphism; allows a good notion of algebraizable logic (\cite{MaMe}).
\qdr\ect

\bct {\em Dense morphism}

$\bullet$ \ $f : l \ra l' \in \Lf$ is dense iff\
$\all \varphi'_n \in F(\Sigma')[n]$ \ $\exists \varphi_n \in F(\Sigma)[n]$ such that\
 $\varphi'_n \dashv'\vdash \check{f}(\varphi_n)$.

$\bullet$ \ $f : l \ra l'$ is  a $\cL$-epimorphism (= surjective at each level $n \in \N$), thus it is a  dense $\cL$-morphism.

$\bullet$ \ $l' \in \Lc$ $\Ra$\
 $f$ is dense iff $\all c'_n \in \Sigma'_n$ \ $\exists \varphi_n \in F(\Sigma)[n]$ such that\
  $c'_n(x_0, \ldots, x_{n-1}) \dashv'\vdash \check{f}(\varphi_n)$.
\qdr\ect

\bct {\em $\Qc$-isos:} \ For $h \in \Lc(l,l')$, are equivalent:

$\bullet$ \ $[h] \in \Qc(l,l')$ is $\Qc$-isomorphism;

$\bullet$ \ $h$ is a dense morphism and $h$ is a conservative translation (i.e., $\Gamma \vdash \psi \ \Lra \ \check{h}[\Gamma] \vdash' \check{h}(\psi)$, for all $\Gamma \cup\{\psi\} \sub F(\Sigma)$).
\qdr\ect

\bct {\em Quotient categories of (Lindenbaum) algebraizable logics} \\
$QLind(\cA_f) \hookr Q\cA_f$:\\
$\bullet$ \  closed under directed colimits\\
$\bullet$ \ reflective subcategory\\
$\bullet$ \ both have non-empty colimits 
\qdr \ect

\begin{picture}(180,70)
\setlength{\unitlength}{.4\unitlength} \thicklines 
\put(-2,20){$Lind(\cA_f)$}
\put(32,162){$\cA_f$} \put(35,154){\vector(0,-1){115}}
\put(45,39){\vector(0,1){115}}
\put(180,154){\vector(0,-1){115}}
\put(190,39){\vector(0,1){115}}
\put(325,154){\vector(0,-1){115}}
\put(335,39){\vector(0,1){115}}
 \put(65,168){\vector(1,0){100}}
  \put(105,25){\vector(1,0){60}}
  \put(177,162){$\Lf$} 
  \put(15,90){$L$} \put(47,90){$j$}
\put(100,175){\small $incl$} 
\put(120,3){\small $incl$} 
\put(177,20){$\Lc$} \put(320,20){$Q\Lc$} \put(320,162){$Q\Lf$}\put(195,90){$i$}
\put(165,90){$c$} \put(340,90){$\bar{i}$} \put(310,90){$\bar{c}$} \put(250,3){$q^c$} \put(250,180){$q$}
\put(210,168){\vector(1,0){100}} \put(210,25){\vector(1,0){100}}
\end{picture}

\subsection{Algebraizable Logics and Categories}

\bct Recall that in the theory of Blok-Pigozzi, to each algebraizable logic $a = (\Sigma, \vdash)$  is canonically associated a {\em unique} quasivariety $QV(a)$ in the same signature $\Theta$ (its "algebraic codification"). 
\qdr\ect

\ble \label{adjQV-le} The inclusion functor has a left adjoint $(L,I): QV \rightleftarrows \alpha-Str$: given by $M \mapsto M/\theta_M$ where $\theta_M$ is the least $\Sigma$-congruence in $M$ such that $M/\theta_M \in QV$. Moreover, the unity of the adjunction $(L,I)$ has components $(q_M)_{M \in \Sigma-Str}$, where $q_M : M \thra M/\theta_M$ is the quotient homomorphism.
\qdr\ele

\bre \label{adjQV-re} The (forgetful) functor $(QV \overset{I}\to\ra \Sigma-Str \overset{U}\to\ra Set)$ has the (free) functor $(Set \overset{F}\to\ra \Sigma-Str \overset{L}\to\ra QV)$, $Y \mapsto F(Y)/\theta_{F(Y)}$, as left adjoint. Moreover, if $\sigma_Y : Y \ra U\circ F(Y)$ is the $Y$-component of the  unity of the adjunction $(F,U)$, then $(Y \overset{t_Y}\to\ra UILF(Y))  \ : = \ (Y \overset{\sigma_Y}\to\ra UF(Y)  \overset{U(q_{F(Y)})}\to\ra U I L F(Y))$ is the $Y$-component of the adjunction $(L\circ F, U \circ I)$.  
\qdr\ere

\bte \ Let $h \in \cA_f(a,a')$, then the induced functor $h^\star : \Sigma'-Str \ra \Sigma-Str$ \
($M'  \overset{}\to\mapsto (M')^h$), "commutes over $Set$" (i.e., $U \circ h^\star = U'$) and has the following additional properties:\\
(a) it has restriction $h^\star\!\!\!\rest : QV(a') \ra QV(a)$ \ 
(i.e. $I \circ h^\star\!\!\!\rest = h^\star \circ I'$);\\
(b) there is a natural epimorphism $\tilde{h}:  L \circ h^\star \thra h^\star\!\!\!\rest \circ L'$, that restricts to $ L \circ h^\star\circ I' = h^\star\!\!\!\rest \circ L' \circ I'$  
\qdr\ete

 \begin{picture}(120,70)
\setlength{\unitlength}{.4\unitlength} \thicklines 
\put(-13,20){\small $\Sigma\!\!-\!\!str$}
\put(-8,162){\small $QVa$} 
\put(160,162){\small $QVa'$} 
\put(40,90){\small $L$}
\put(05,90){\small $I$} 
\put(100,175){\small $h^\star\!\!\rest$} 
\put(155,20){\small $\Sigma'\!\!-\!\!str$} 
\put(100,90){$\nearrow$} 
\put(105,70){\small $\bar{h}$} 
\put(155,90){\small $L'$} 
\put(200,90){\small $I'$} 
\put(100,0){\small $h^\star$}
\put(25,150){\vector(0,-1){100}}
\put(195,150){\vector(0,-1){100}} 
\put(35,50){\vector(0,1){100}}
\put(185,50){\vector(0,1){100}} 
\put(150,167){\vector(-1,0){80}} 
\put(150,25){\vector(-1,0){80}}
\put(35,-50){\small $U$}
\put(165,-50){\small $U'$} 
\put(75,-40){\small $F$}
\put(125,-40){\small $F'$} 
\put(95,-130){\small $Set$} 
\put(102,-105){\vector(-1,2){55}} 
\put(118,-105){\vector(1,2){55}}
\put(35,05){\vector(1,-2){55}} 
\put(185,05){\vector(-1,-2){55}}
\qbezier(20,150)(-100,15)(70,-115) 
\put(65,-110){\vector(1,-1){05}}
\put(155,-110){\vector(-1,-1){05}}
\qbezier(200,150)(320,15)(150,-115) 
\end{picture}
 
\vspace{2cm}

{\bf Good representation theory of $Lind(\cA_f)$}

\bpr  Let  $g_0, g_1 : l \ra a \in \cL_f$, with $a \in Lind(\cA_f)$.

(a) \  $g_0$ is dense $\Ra$ $g^\star\!\!\rest : QV(a) \ra \Sigma-str$ is full, faitful and injective on objects.

(b) \ $[g_0] = [g_1] \in \Qf$ $\Ra$ $g_0^\star\!\!\rest = g_1^\star\!\!\rest : QV(a) \ra \Sigma-str$.
\qdr\epr

\bpr \label{lindalg-ct} Let $a = (\Sigma, \vdash)$ be a Lindenbaum algebraizable then:\\
(a) \ $F(\Sigma)/\Delta  = F(\Sigma)/\!(\dashv \ \vdash)$ is a $\Sigma$-structure.\\
(b) \ $F(\Sigma)/\Delta \in QV(a)$.\\
(c) \ $F(\Sigma)/\Delta$ is the free $QV(a)$-object over the set $X = \{x_0, \ldots, x_n, \ldots\}$.
 \qdr
\epr

\bpr \label{lindalgiso-pr} Let $a$ and $a'$ be Lindenbaum algebraizable logics. If  $a \overset{[h']}\to{\underset{[h]}\to\rightleftarrows} a'$ is a pair of inverse $QLind(\cA_f)$-isomorphisms (i.e., are $\Qc$-isomorphisms that preserve algebraizing pair) then:\ $QV(a) \overset{{h'}^\star\!\!\rest}\to{\underset{{h}^\star\!\!\rest}\to\leftrightarrows} QV(a')$ is an isomorphism of categories. 
\qdr\epr

\ble  Let $\Sigma, \Sigma' \in Obj(\Sf)$. Consider $H: \Sigma'-Str \ra \Sigma-Str$ a functor that "commutes over $Set$" (i.e. $U \circ H = U'$) and, for each set $Y$, let $\eta_H(Y) : F(Y) \ra H(F'(Y))$ be the unique $\Sigma$-morphism such that $(Y \overset{\sigma_Y}\to\ra UF(Y)  \overset{U(\eta_Y)}\to\ra UHF'(Y)) \ = \ (Y \overset{\sigma'_Y}\to\ra U'F'(Y))$. Then:\\ 
(a)\ For each set $Y$ and each $\psi \in F(Y)$, $Var(\eta_H(Y)(\psi)) \ \sub \ Var(\psi)$;\\  
(b)\ $(\eta_H(Y))_{Y \in Set}$ is a natural transformation $\eta_{H}: F \ra H\circ F'$;\\
(c)\  If $H$ is an isomorphism of categories, then $\eta_H(Y)$ "preserves variables" (i.e., $\all \psi \in F(Y)$, $Var(\eta_H(Y)(\psi)) \ = \ Var(\psi)$) and $H$ preserves (strictly) products and substructures.\\
(d)\ For each  $n \in \N$, let $X_n := \{x_0, \cdots, x_{n-1}\} \sub X$, if $\eta_H(X_n)$ "preserves variables", then the mapping $c_n \in \Sigma_n \ \mapsto \ \eta_H(X_n)(c_n(x_0, \cdots, x_{n-1})) \in F'(X_n)$ determines a $\Sf$-morphism $m_H : \Sigma \ra \Sigma'$.
\qdr\ele

\bct \\
(a) Let $\Sigma, \Sigma' \in Obj(\Sf)$. Let $H : \Sigma'-Str \ra \Sigma-Str$ be a  "signature" functor, i.e. a functor satisfying $(s1), (s2), (s3)$: \\
\hfl $(s1)$ $H$ "commutes over $Set$" \hfl \hfl $(s2)$ $\eta_H$ "preserves variables" \hfl \\
\hfl $(s3)$ $H$ preserves (strictly) products and substructures.\hfl\\
(b) Denote $\Sf^\dagger$ the subcategory of the category of diagrams (i.e., the category whose objects are categories and the arrows are change of base morphisms (i.e., some pairs (functors, natural transformations)), given by all the categories $\Sigma-str$  and morphisms $(H, \eta_H)$ where $H$ is a signature functor.\\
(c) Let $a, a' \in Obj(Lind(\cA_f))$.  Let $H : \Sigma'-Str \ra \Sigma-Str$ be a  "Lindenbaum" functor, i.e. a signature functor also satisfying $(l1), (l2), (l3)$: \\
\hfl $(l1)$ $H$ has a (unique) restriction to the quasivarieties $H\rest : QV(a') \ra QV(a)$ \hfl \\
\hfl $(l2)$ $\check{m}_{H}(\Delta)\dashv ' \vdash\Delta'$ \hfl \hfl $(l3)$ $\check{m}_{H}(\delta)= \check{m}_H(\varepsilon)\vDash_{QV(a')}\vDash \delta'= \varepsilon'$ \hfl.\\
(d) Let $a, a' \in Obj(Lind(\cA_f))$ and $H : \Sigma'-Str \ra \Sigma-Str$ be a  "Lindenbaum" functor. For each set $Y$, let $\bar{\eta}_H(Y) : LF(Y) \ra H\rest(L'F'(Y))$ be the unique $QV(a)$-morphism such that\ $(Y \overset{t_Y}\to\ra UILF(Y)  \overset{UI(\bar{\eta}_Y)}\to\ra UIH\rest L' F'(Y)) \ =$ \\ $(Y \overset{t'_Y}\to\ra U'I'L'F'(Y))$. Then  $(\bar{\eta}_H(Y))_{Y \in Set}$ is a natural transformation $\bar{\eta}_{H}: L \circ F \ra H\rest \circ L'\circ F'$.\\ 
(d) Denote $Lind(\cA_f)^\dagger$ the subcategory of the category of diagrams, given by all the subcategories $QV(a) \hookr \Sigma-str$  and morphisms $(H\rest, \bar{\eta}_H)$ where $H$ is a Lindenbaum functor.
\qdr\ect

\bte 

(a) The categories $\Sf$ and $\Sf^\dagger$ are anti-isomorphic. More precisely, given $\Sigma, \Sigma' \in \Sf$, the mappings $h \in \Sf(\Sigma,\Sigma') \ \mapsto\ (h^\star, \eta_{h^\star}) \in \Sf(\Sigma'-str,\Sigma-str)^\dagger$ and $(H, \eta_{H}) \in \Sf(\Sigma'-str,\Sigma-str)^\dagger\ \mapsto\ m_H \in \Sf(\Sigma,\Sigma')$   are inverse bijections.

(b) The pair of inverse anti-isomorphisms above restricts to a pair of inverse anti-isomorphisms between the categories $Lind(\cA_f)$ and $Lind(\cA_f)^\dagger$. 

Moreover,  the inverse isomorphisms establish a correspondence:

(c) If $h \in Lind(\cA_f)(a,a')$ and $H \in Lind(\cA_f)^\dagger$ are in correspondence, then also are in correspondence the equivalence class $\{h' \in Lind(\cA_f)(a,a'): [h]= [h']  \in QLind(\cA_f)(a,a')\}$ and   the   equivalence class $\{H' \in Lind(\cA_f)^\dagger : H'\rest = H\rest$, $\bar{\eta}_{H'} = \bar{\eta}_{H}\}$. 

(d)  If $h \in Lind(\cA_f)(a,a')$ and $H \in Lind(\cA_f)^\dagger$ are in correspondence, then $[h]$ is a $QLind(\cA_f)$-isomorphism (i.e., $h$ is an equipolence of logics) \hem $\Lra$ \hem $(H\!\!\rest, \bar{\eta}_H)$ is an   isomorphism of change of bases.
\qdr\ete






\section{Representation Theory of Logics}

\subsection{General Logics and Categories}



Let $U : \cL_s \ra Lind(\cA_s)$ denote the forgetful functor. 

$\bullet$ {\bf Objects}: To each logic $l=(\Sigma,\vdash)$, are associated two pairs (left and right) of data:

{\bf (I)} two comma categories (over $Lind(\cA_s)$): \\
\bu( $l \ra U$), the "left algebrizable spectrum of $l$" (analysis process); \\
\bu ($U \ra l$),  the "right algebrizable spectrum of $l$"  (synthesis process).

 
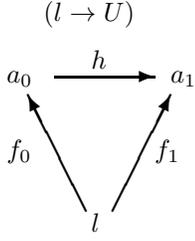
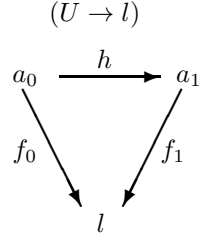
\begin{picture}(120,100)
\thicklines \setlength{\unitlength}{.4\unitlength} 
\put(40,220){$(l \ra U)$}
\put(05,162){$a_0$}
\put(160,162){$a_1$} 
\put(05,90){$f_0$}
\put(145,90){$f_1$} 
\put(85,20){$l$} 
\put(85,175){$h$}
\put(50,167){\vector(1,0){95}} 
\put(80,40){\vector(-1,2){55}} 
\put(105,40){\vector(1,2){55}}

\put(640,220){$(U \ra l)$}
\put(605,162){$a_0$}
\put(760,162){$a_1$} 
\put(605,90){$f_0$}
\put(745,90){$f_1$} 
\put(685,20){$l$} 
\put(685,175){$h$}
\put(650,167){\vector(1,0){95}} 
\put(615,155){\vector(1,-2){55}} 
\put(765,155){\vector(-1,-2){55}}
\end{picture}


 {\bf (II)} two diagrams (left and right "representation diagram"):\\ 
 \bu $l$-$Mod$ $\lrsa$ ($l \ra U, I)$; \\
 \bu $Mod$-$l$ $\lrsa$ ($U \ra l, L)$.

 $l$-$Mod$ : $(l \ra U)^{op} \ \lra \ (\Sigma-str \la CAT)$ \\
 $(a_0, f_0) \hem \mapsto \hem (\Sigma-str \overset{f_0^\star\!\!\rest I_0}\to\la QV(a_0))$\\
 $(a_0,f_0) \overset{h}\to\mapsto$ $ (a_1,f_1) \hem \mapsto$\
 $ \hem ((QV(a_1), {f_1^\star\!\!\rest I_1}) \overset{h^\star\!\!\rest}\to\ra (QV(a_0),{f_0^\star\!\!\rest I_0}))$

 
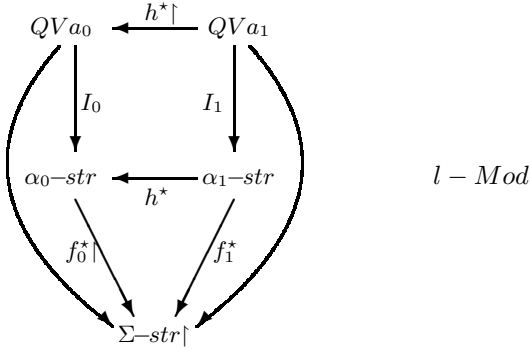
\begin{picture}(120,80)
\setlength{\unitlength}{.4\unitlength} \thicklines 
\put(-13,20){\small $\alpha_0\!\!-\!\!str$}
\put(-8,162){\small $QVa_0$} 
\put(160,162){\small $QVa_1$} 
\put(40,90){\small $I_0$} 
\put(100,175){\small $h^\star\!\!\rest$} 
\put(155,20){\small $\alpha_1\!\!-\!\!str$} 
\put(365,20){ $l-Mod$} 
\put(155,90){\small $I_1$} 
\put(100,0){\small $h^\star$}
\put(35,150){\vector(0,-1){100}}
\put(185,150){\vector(0,-1){100}} 
\put(150,167){\vector(-1,0){80}} 
\put(150,25){\vector(-1,0){80}}
\put(25,-50){\small $f^\star_0\!\!\rest$}
\put(165,-50){\small $f^\star_1$} 
\put(75,-130){\small $\Sigma\!\!-\!\!str\!\!\rest$} 
\put(35,05){\vector(1,-2){55}} 
\put(185,05){\vector(-1,-2){55}}
\qbezier(20,150)(-100,15)(70,-115) 
\put(65,-110){\vector(1,-1){05}}
\put(155,-110){\vector(-1,-1){05}}
\qbezier(200,150)(320,15)(150,-115) 
\end{picture}
 
 \vspace{3cm}

 $Mod$-$l$ : $(U \ra l) \ \lra \ (2-CAT \la \Sigma-str )$ \\
 $(a_0, f_0) \hem \mapsto \hem (QV(a_0) \overset{L_0 f_0^\star\!\!\rest}\to\la  \Sigma-str)$\\
  $(a_0,f_0) \overset{h}\to\mapsto$
  $ (a_1,f_1) \hem \mapsto$\
  $ \hem ((QV(a_0), {L_0 f_0^\star\!\!\rest}) \overset{(h^\star\!\!\rest, \tilde{h})}\to\ra (QV(a_1),{L_1 f_1^\star\!\!\rest}))$

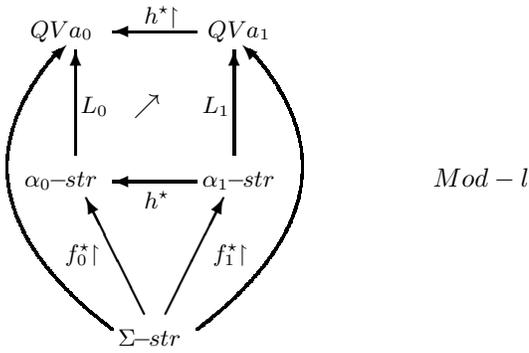
 \begin{picture}(120,80)
\setlength{\unitlength}{.4\unitlength} \thicklines 
\put(-13,20){\small $\alpha_0\!\!-\!\!str$}
\put(-8,162){\small $QVa_0$} 
\put(160,162){\small $QVa_1$} 
\put(40,90){\small $L_0$} 
\put(100,175){\small $h^\star\!\!\rest$} 
\put(155,20){\small $\alpha_1\!\!-\!\!str$} 
\put(365,20){ $Mod-l$} 
\put(155,90){\small $L_1$} 
\put(100,0){\small $h^\star$}
\put(35,50){\vector(0,1){100}}
\put(185,50){\vector(0,1){100}} 
\put(150,167){\vector(-1,0){80}} 
\put(150,25){\vector(-1,0){80}}
\put(25,-50){\small $f^\star_0\!\!\rest$}
\put(165,-50){\small $f^\star_1\!\!\rest$} 
\put(75,-130){\small $\Sigma\!\!-\!\!str$} 
\put(100,-100){\vector(-1,2){55}} 
\put(120,-100){\vector(1,2){55}}
\qbezier(20,150)(-100,15)(70,-115) 
\put(22,150){\vector(1,1){05}}
\put(198,150){\vector(-1,1){05}}
\qbezier(200,150)(320,15)(150,-115) 
\put(90,90){$\nearrow$}
\end{picture}
 
 \vspace{2cm}
 
$L_0f_0^\star \overset{\tilde{h}f_1^\star}\to\Ra h^\star\!\!\rest L_1f_1^\star$\hem
where\hem
$L_0h^\star \overset{\tilde{h}}\to\Ra h^\star\!\!\rest L_1$

Since:\hem
$(L_0h^\star g^\star \overset{\tilde{h}g^\star}\to\Ra h^\star\rest L_1g^\star
 \overset{h^\star\!\!\rest\tilde{g}}\to\Ra h^\star\!\!\rest  g^\star\!\!\rest L_2)$\ = \
$L_0 (gh)^\star \overset{\widetilde{gh}}\to\Ra (gh)^\star\!\!\rest L_2$\\
\hp then: \hem
$(L_0 f_0^\star \overset{\tilde{h}f_1^\star}\to\Ra h^\star\!\!\rest L_1f_1^\star
 \overset{h^\star\!\!\rest(\tilde{g}f_2^\star)}\to\Ra h^\star\!\!\rest  g^\star\!\!\rest L_2f_2^\star)$ \ =\ 
$L_0 (gh)^\star \overset{\widetilde{gh}f_2^\star}\to\Ra (gh)^\star\!\!\rest L_2f_2^\star$

$\bullet$ {\bf Arrows}:  Morphisms between logics $t : l \ra l'$ induce two pairs (left and right) of data: \


{\bf (I)} a left/right "spectral" functor:  \\ 
$(l \ra U) \ \overset{- \circ t}\to\la  \ (l' \ra U)$; \\
$( U \ra l) \ \overset{t \circ -}\to\ra  \ (U \ra l')$.


{\bf (II)} a left/right "representation diagram" morphism: \\
$(l$-$Mod) \ \overset{t^\star \circ -}\to\la  \ (l'$-$Mod)$; \\
$(Mod$-$l) \ \overset{- \circ t^\star}\to\ra  \ (Mod$-$l')$. 

 \begin{picture}(120,80)
\thicklines \setlength{\unitlength}{.4\unitlength} 
\put(-25,162){\small $QVa_0$}
\put(150,162){\small $QVa_1$} 
\put(40,130){\small $f^\star_0I_0$}
\put(95,110){\small $f^\star_1I_1$} 
\put(60,20){\small $\Sigma'\!\!-\!\!str$}
\put(350,60){$(l$-$Mod) \ \overset{t^\star \circ -}\to\la  \ (l'$-$Mod)$}
\put(60,-80){\small $\Sigma\!\!-\!\!str$} 
\put(85,175){\small $h^\star\!\!\rest$}
\put(70,-25){\small $t^\star$} 
\put(140,167){\vector(-1,0){105}} 
\put(15,155){\vector(1,-2){55}} 
\put(165,155){\vector(-1,-2){55}}
\put(90,10){\vector(0,-1){60}}
\put(5,150){\vector(1,-3){68}} 
\put(175,150){\vector(-1,-3){68}}
\end{picture}

\vspace{1cm}

\bct {\em The category of all left modules: {\bf LM}}

\und{objects}: a left module for a logic $l$, i.e. the functor $left(l) : (l \ra U)^{op} \ \overset{l-Mod}\to\lra \ (\Sigma-str \la CAT) \ \hookr (CAT \la CAT)$

\und{arrows}: a pair $(B, \tau) : left(l') \ra left(l)$, where $B :  (l' \ra U) \ra (l \ra U)$ is a "change of bases" functor, and $\tau : left(l') \Ra left(l) \circ B$ is a natural transformation with the additional compatibility condition:\\
 for each  $(a'_0, f'_0), (a'_1, f'_1) \in (l' \ra U) $ \ $\Ra$ \ $Proj(\tau_{(a'_0, f'_0)}) = Proj(\tau_{(a'_1, f'_1)}) : \Sigma'-str \ra \Sigma-str$
\qdr\ect

\begin{picture}(200,130)
\setlength{\unitlength}{.4\unitlength} \thicklines 
\put(-45,304){\small $QVcod(f'_1)$} 
\put(-05,232){\small ${f'_1}^\star$} 
\put(155,304){\small  $QVcod(B(f'_1))$} 
\put(200,232){\small $B(f'_1)^\star$} 
\put(90,227){\small $\tau(f'_1)$}
\put(100,247){$\Ra$}
\put(25,292){\vector(0,-1){100}}
\put(195,292){\vector(0,-1){100}} 
\put(75,309){\vector(1,0){80}} 

\put(-13,162){\small $\Sigma'\!\!-\!\!str$}
\put(-45,20){\small $QVcod(f'_0)$} 
\put(170,162){\small $\Sigma\!\!-\!\!str$} 
\put(-5,90){\small ${f'_0}^\star$}
\put(155,20){\small  $QVcod(B(f'_0))$} 
\put(200,90){\small $B(f'_0)^\star$} 
\put(90,105){\small $\tau(f'_0)$}
\put(100,85){$\Ra$}
\put(25,50){\vector(0,1){100}}
\put(195,50){\vector(0,1){100}} 
\put(55,167){\vector(1,0){110}} 
\put(75,25){\vector(1,0){80}}
\end{picture}

\bpr \\
\hp  (a) $t: l \overset{\cong}\to\ra l' \in \Ls$ \ $\Ra$ \
$(-\circ t, (t^\star, id)) : left(l') \overset{\cong}\to\ra left(l) \in\ ${\bf LM}. 

(b) $can_l: l \ra l^{(c)}$ induces a {\bf LM}-isomorphism: $left(l^{(c)}) \overset{\cong}\to\ra left(l)$
\qdr \epr


\subsection{Morita equivalence of logics and variants}

$\bullet$ \ (left/right) Morita equivalence of rings: an equivalence relation coarser than isomorphism \\
({\small  Ex.: For rings, $R \equiv Mat_{n\times n}(R)$})

\bdf The logics $l$ and $l'$ are left Morita equivalent when:

(a) Let $S$ be a full subcategory of $(l \ra U)$. $S$ is called {\em generic} if $S^{op} \hookr (l \ra U)^{op} \overset{l-Mod}\to\ra (CAT \ra \Sigma-str)$ is "relatively cofinal" (in the image...)

(b) The logics $l$ and $l'$ are left Morita equivalent when there are:\\
$\bullet$ generic subcategories $S  \hookr (l \ra U)$ and $S' \hookr (l' \ra U)$;\\
$\bullet$ functors $B : S' \ra S$ and $B': S \ra S'$;\\
$\bullet$ "natural comparations": $(T,\tau)$ and $(T',\tau')$ \\
\hp $\ast$ $T: \Sigma'-str \ra \Sigma-str$ is a functor and for each $(a',f') \in S'$, $\tau_{f'} : QVcod(f') \overset{\cong}\to\ra QVcod(B(f'))$ is a isomorphism of categories such that $B(f')^{\star}  \circ \tau_{f'} = T \circ f'^\star$ and for each $(a_0',f_0') \overset{g}\to\ra (a_1',f_1') \in S'$\\
\hp $\ast$ analogous conditions for $(T',\tau')$ 

\begin{picture}(120,80)
\setlength{\unitlength}{.4\unitlength} \thicklines 
\put(-13,162){\small $\Sigma'\!\!-\!\!str$}
\put(-25,20){\small $QVcod(f')$} 
\put(180,162){\small $\Sigma\!\!-\!\!str$} 
\put(-10,90){\small $f'^\star$}
\put(100,175){\small $T$} 
\put(165,20){\small  $QVcod(B(f'))$} 
\put(200,90){\small $B(f')^\star$} 
\put(110,05){\small $\tau_{f'}$}
\put(20,50){\vector(0,1){100}}
\put(195,50){\vector(0,1){100}} 
\put(65,167){\vector(1,0){100}} 
\put(95,25){\vector(1,0){60}}

\end{picture}

\begin{picture}(200,120)
\setlength{\unitlength}{.35\unitlength} \thicklines 
\put(-25,304){\small $QVcod(f'_1)$} 
\put(95,232){\small ${f'_1}^\star$} 
\put(425,304){\small  $QVcod(B(f'_1))$} 
\put(330,232){\small $B(f'_1)^\star$} 
\put(220,322){$\tau(f'_1)$}
\put(240,282){$\cong$}
\put(25,292){\vector(1,-1){110}}
\put(465,292){\vector(-1,-1){110}} 
\put(475,295){\vector(0,-1){250}}
\put(105,309){\vector(1,0){310}} 

\put(405,162){\small ${B(g')}^\star$}
\put(20,162){\small ${g'}^\star$}
\put(310,162){\small $\Sigma\!\!-\!\!str$}
\put(-25,20){\small $QVcod(f'_0)$} 
\put(250,172){\small $T$}
\put(145,162){\small $\Sigma'\!\!-\!\!str$} 
\put(95,90){\small ${f'_0}^\star$}
\put(425,20){\small  $QVcod(B(f'_0))$} 
\put(330,90){\small $B(f'_0)^\star$} 
\put(220,00){$\tau(f'_0)$}
\put(240,30){$\cong$}
\put(25,45){\vector(1,1){110}}
\put(15,295){\vector(0,-1){250}}
\put(465,45){\vector(-1,1){110}} 
\put(215,167){\vector(1,0){90}} 
\put(105,25){\vector(1,0){310}}
\end{picture}

\edf

\bte If $left(l) \cong left(l')$ then  $l$ and $l'$ are left Morita equivalent. In particular:\\
(a) If $l \cong l'$, then $l$ and $l'$ are left Morita equivalent.\\
(b) $l$ and $l^{(c)}$ are left Morita equivalent.
\qdr\ete





\bte If $l$ and $l'$ are equipollent, then they are left Morita equivalent.
\qdr\ete 

\ble Let $l \overset{[t]}\to\ra l'$ be a  $\Qf$-isomorphism. Then:

(a) \ For each Lindenbaum algebraizale logic $a'= (\alpha', \vdash')$ and any dense $\Ls$-morphism $f': l' \ra a'$, consider:\\
$\bullet$ \ $\alpha = (\Sigma_n/\equiv_n)_{n \in \N}$, where $c_n \equiv_n d_n$ iff $\hat{f}(\check{t}(c_n)) = \hat{f}(\check{t}(d_n)$; \\
$\bullet$ \ $f = quoc$;  \\
$\bullet$ \  
$h$ be the unique $\Sf$-morphim such that
$h \bullet f = f' \bullet t$ in $\Sf$ (thus \ $\check{h} \circ \hat{f} = \hat{f'} \circ \check{t}$)\\
$\bullet$ \ $a = (\alpha, h^\star(\vdash_{a'}))$;\\
Then:\\
$\bullet$ \ $a$ is a Lindenbaum algebraizable logic;\\
$\bullet$ \ $h : a \ra a'$ is a $\Lc$-morphism preserves algebraizable pair and is an weak equivalence i.e. it induces a  $\Qc$-isomorphism;   \\
$\bullet$ \ $f : l \ra a$ is  a $\Ls$-epimorphism (= surjective at each level $n \in \N$), thus it is a  dense $\Ls$-morphism.

(b)  If $a_0'= (\alpha_0', \vdash_0')$, $a_1'= (\alpha_1', \vdash_1')$ are Lindenbaum algebraizale logic;  $f'_0: l' \ra a_0'$, $f'_1: l' \ra a_1'$ are dense $\Ls$-morphisms and  $g': a'_0 \ra a_1'$ is a $\Ls$-morphism such that $g'f'_0 = f'_1$, then: \\
$\bullet$ \ there is a unique $\Ss$-morphism $g : \alpha_0 \ra \alpha_1$ such that $gf_0 = f_1$;
 Moreover:\\
$\bullet$ \ $h_1 g = g' h_0$;\\
$\bullet$ \ $g$ is a $\cA_s$-morphism.
\qdr\ele

\bdf The logics $l$ and $l'$ are left-stably Morita equivalent when:\\
$\bullet$ \ there are functors: \ $\Sigma'-str \overset{F}\to{\underset{F'}\to\rlas} \Sigma-str$;\\
$\bullet$ \ there are functors:\
 $colim_{ f' \in (l'\ra U)} \ QVcod(f') \overset{E}\to{\underset{E'}\to\rlas} colim_{f \in (l\ra U)}\  QVcod(f)$;\\
such that:\\
$\bullet$ \ $E$ and $E'$ are quasi-inverse equivalence functors;\\
$\bullet$ \ the diagram below commutes:



\begin{picture}(120,80)
\setlength{\unitlength}{.4\unitlength} \thicklines 
\put(-03,20){\small $\Sigma'\!\!-\!\!str$}
\put(02,162){\small $\Sigma\!\!-\!\!str$} 
\put(350,162){\small $colim\ QV(cod(f))$} 
\put(-15,90){\small $F$}
\put(55,90){\small $F'$} 
\put(200,175){\small $can$} 
\put(350,20){\small $colim\ QV(cod(f'))$} 
\put(450,90){\small $E$} 
\put(340,90){\small $E'$} 
\put(200,0){\small $can'$}
\put(395,90){\small $\simeq$}
\put(45,150){\vector(0,-1){100}}
\put(370,150){\vector(0,-1){100}} 
\put(15,50){\vector(0,1){100}}
\put(430,50){\vector(0,1){100}} 
\put(325,167){\vector(-1,0){240}} 
\put(325,25){\vector(-1,0){240}}
\end{picture}
\qdr\edf

\bpr\ $l$ and $l'$ are left Morita equivalent \ $\Ra$ \ $l$ and $l'$ are left-stably Morita equivalent.
\qdr\epr

\bco\ $left(l) \cong left(l')$ \ $\Ra$ \ $l$ and $l'$ are left-stably Morita equivalent.
\qdr\eco

\bpr\ $a \cong a' \in QLind(A_f)$ \ $\Ra$\
$a$ and $a'$ are left stably Morita equivalent.
\qdr\epr

\bco \label{classMor-co} The presentations of classical logics are left stably Morita equivalent.

\begin{picture}(120,80)
\setlength{\unitlength}{.4\unitlength} \thicklines 
\put(-23,20){\small $(\neg',\vee')\!\!-\!\!str$}
\put(-23,162){\small $(\neg,\ra)\!\!-\!\!str$} 
\put(350,162){\small $BA(\neg,\ra) \cong colim\ QV(cod(f)) $} 
\put(-15,90){\small $t^\star$}
\put(55,90){\small $t'^\star$} 
\put(200,175){\small $incl$} 
\put(350,20){\small $BA(\neg',\vee') \cong colim\ QV(cod(f'))$} 
\put(450,90){\small $t^\star\!\!\rest$} 
\put(320,90){\small $t'^\star\!\!\rest$} 
\put(200,0){\small $incl'$}
\put(395,90){\small $\cong$}
\put(45,150){\vector(0,-1){100}}
\put(370,150){\vector(0,-1){100}} 
\put(15,50){\vector(0,1){100}}
\put(430,50){\vector(0,1){100}} 
\put(340,167){\vector(-1,0){220}} 
\put(340,25){\vector(-1,0){220}}
\end{picture}
\qdr \eco

\bpr \label{classintuit-re} Concerning the relations between Classical logics and Intuitionist logics:

(a) \ They are not stably-Morita-equivalent. 

(b) \ But they are only stably-Morita-adjointly related: \\
$L : HA \ra BA$ : $H \mapsto H/F_H$, where 
$F_H = \ll \{ a \leftrightarrow \neg\neg a: a \in H \} \rr$



\begin{picture}(120,90)
\setlength{\unitlength}{.4\unitlength} \thicklines 
\put(-23,20){\small $(\neg,\vee, \wedge, \ra)\!\!-\!\!str$}
\put(-23,162){\small $(\neg,\vee, \wedge, \ra)\!\!-\!\!str$} 
\put(350,162){\small $HA \cong colim\ QV(cod(f)) $} 
\put(-15,90){\small $id$}
\put(55,90){\small $id$} 
\put(240,175){\small $incl$} 
\put(350,20){\small $BA \cong colim\ QV(cod(f'))$} 
\put(440,90){\small $J$} 
\put(340,90){\small $L$} 
\put(240,0){\small $incl'$}
\put(380,90){\small $adj$}
\put(45,150){\vector(0,-1){100}}
\put(370,150){\vector(0,-1){100}} 
\put(15,50){\vector(0,1){100}}
\put(430,50){\vector(0,1){100}} 
\put(340,167){\vector(-1,0){150}} 
\put(340,25){\vector(-1,0){150}}
\end{picture}
\qdr\epr




\section{Final Remarks and Future Works}

$\bullet$ \ Present the adequade definitions  "on the right side"  that allow get basic results analogous to "left side". Note that the considerations  "on left" and "on right" are  distincts conceptually ("left"' is adequate for  analysis of logics; "right" is related to synthesis process) and technically (a 2-categorial aspect is needed "on right").\\
$\bullet$ \ Describe necessary/sufficient conditions for Morita equivalence of logics (and variants).\\
$\bullet$ \ Induce new (functorial) morphisms between logics from the representation diagrams $left(l)$ and $right(l)$.\\
$\bullet$ \ Analise categories of fractions of  categories of logics.\\
$\bullet$ \ Define a general theory of "mathematical invariants" to mesure the degree of distinctions of arbitrary logics and 
develope general methods of calculation of invariants (in some sense).\\
$\bullet$ \ Understand new notions of identity of logics.\\
$\bullet$ \ Describe similar construction on alternative base categories (\cite{MaMe}). ({\small Example: the study of LFIs  by Possible Translations Semantics.})



\end{document}